\documentclass{amsart}
\usepackage{amssymb}
\usepackage{xy}

\theoremstyle{plain}

\newtheorem{thm}{Theorem}[section]

\newtheorem{cor}[thm]{Corollary}
\newtheorem{prop}[thm]{Proposition}
\newtheorem{lem}[thm]{Lemma}

\newtheorem{con}[thm]{Conjecture}

\theoremstyle{definition}
\newtheorem{defn}[thm]{Definition}
\newtheorem{ex}[thm]{Example}

\input xy
\xyoption{matrix}
\xyoption{arrow}

\title[Lazarsfeld-Mukai bundles]{Lazarsfeld-Mukai bundles and applications}
\author[Aprodu]{Marian Aprodu}
\subjclass[2000]{}
\address{Romanian Academy, 
Institute of Mathematics "Simion Stoilow" 
P.O. Box 1-764, RO 014700, 
Bucharest, Romania}
\email{marian.aprodu@imar.ro}
\address{\c Scoala Normal\u a Superioar\u a Bucure\c sti, 
Calea Grivi\c tei 21, RO-010702, Bucharest, Romania}
\date{\today}                                           
\thanks{This work was partly supported by 
the grant "Vector bundle techniques in the geometry of complex varieties", 
PN-II-ID-PCE-2011-3-0288, contract no. 132/05.10.2011. I am 
grateful to S. Druel, G. Farkas and A. Ortega for useful discussions
on this subject.}

\begin{document}

\begin{abstract}
 We survey the development of the notion of Lazarsfeld-Mukai bundles
together with various applications, from the classification of Mukai manifolds
to Brill-Noether theory  and syzygies of $K3$ sections. 
To see these techniques at work, we present a short
proof of a result of M. Reid on the existence of elliptic pencils.
\end{abstract}
\maketitle

\section*{Introduction}

Lazarsfeld--Mukai bundles appeared naturally in connection with two 
completely different important problems in algebraic geometry from the 
1980s. The first problem, solved by Lazarsfeld, was to find explicit examples 
of smooth curves which are generic in the sense of Brill-Noether-Petri
\cite{LazarsfeldJDG}. The second problem was the classification
of prime Fano manifolds of coindex 3 \cite{MukaiPNAS}. More recently,
Lazarsfeld--Mukai bundles have found
applications to syzygies and higher-rank Brill--Noether theory.

The common
feature of all these research topics is the central role played by
$K3$ surfaces and their hyperplane sections. For the Brill--Noether--Petri genericity, Lazarsfeld
proves that a general curve in a linear system that generates
the Picard group of a $K3$ surface satisfies this condition.
For the classification of prime Fano manifolds of coindex 3,
after having proved the existence of smooth fundamental divisors,
one uses the geometry of a two-dimensional linear section which is a very general
$K3$ surface.

The idea behind this definition is that the Brill--Noether theory
of smooth curves~on a $K3$ surface, also called \textit{$K3$ sections},
is governed by higher-rank vector bundles on the surface. To be
more precise, consider $S$ a $K3$ surface (considered always to be
smooth, complex, projective), $C$ a smooth curve on $S$ of genus $\ge 2$,
and $|A|$ a base-point-free pencil on $C$.
If we attempt to lift the linear system
$|A|$ to the surface $S$, in~most cases, we will fail.
For instance, $|A|$ cannot
lift to a pencil on $S$ if $C$ generates
$\mathrm{Pic}(S)$ or if $S$ does not contain any elliptic curve at all.
However, interpreting a general divisor in $|A|$
as a zero-dimensional subscheme of $S$, it is
natural to try and find a rank-two bundle $E$ on $S$ and
a global section of $E$ whose scheme of zeros coincides with
the divisor in question. Varying the divisor, one~should exhibit
in fact a two-dimensional space of global sections of $E$.
The effective construction of $E$ is realized through
elementary modifications, see Sect.\,\ref{section: Def}, and this
is precisely a Lazarsfeld--Mukai bundle of rank two.
The passage to higher ranks is natural, if we start
with a complete, higher-dimensional, base-point-free linear
system on $C$. At the end, we obtain vector bundles
with unusually high number of global sections, which
provide us with a rich geometric environment.

\medskip

The structure of this chapter is as follows.
In the first section, we recall the definition of Lazarsfeld--Mukai
bundles and its first properties.
We note equivalent conditions for a bundle to be Lazarsfeld--Mukai in
Sect.\,\ref{subsection: Definition LM},
and we discuss simplicity in the rank-two case in
Sect.\,\ref{subsection: Simple LM}.
The relation with the Petri conjecture and the
classification of Mukai manifolds, the original
motivating problems for the~definition, are  considered in
Sects.\,\ref{subsection: BNP} and \ref{subsection: Mukai manifolds}, respectively.
In Sect.\,\ref{section: Constancy}
we treat the problem of constancy of invariants in
a given linear system. For small gonalities, Saint-Donat
and Reid proved that minimal pencils on $K3$ sections
are induced from elliptic pencils on the $K3$ surface;
we present a short proof using Lazarsfeld--Mukai bundles
in Sect.\,\ref{subsection Gonality I}.
Harris and Mumford conjectured that the gonality
should always be constant. We discuss the evolution
of this conjecture, from Donagi--Morrison's counterexample,
Sect.\,\ref{subsection Gonality I},
to Green--Lazarsfeld's reformulation in terms
of Clifford index, Sect.\,\ref{subsection: Cliff} and to Ciliberto--Pareschi's results on the subject, Sect.\,\ref{subsection: Gonality II}.
The works around this problem emphasized the importance
of parameter spaces of Lazarsfeld--Mukai bundles.
We conclude the section with a discussion of
dimension calculations of these spaces, Sect.\,\ref{subsection: Parameter}, which are
applied afterwards to Green's conjecture.
Sect.\,\ref{section: Green} is devoted to Koszul cohomology
and notably to Green's conjecture for $K3$ sections.
After recalling the definition and the motivations that
led to the definition, we discuss the statement of Green's
conjecture, and we sketch the proof for $K3$ sections.
Voisin's approach using punctual Hilbert schemes, which is
an essential ingredient, is examined in Sect.\,\ref{subsection: Hilbert}. Lazarsfeld--Mukai bundles
are fundamental objects in this topic, and their
role is outlined in Sect.\,\ref{subsection: Role of LM}.
The final step in the solution of Green's conjecture for
$K3$ sections is tackled in Sect.\,\ref{subsection: Green for K3}.
We conclude this chapter with a short discussion on
Farkas--Ortega's new applications of Lazarsfeld--Mukai bundles to
Mercat's conjecture (which belongs to the rapidly developing
higher-dimensional Brill--Noether theory),
Sect.\,\ref{section: Higher BN}.

\medskip

\noindent{\it Notation.} The additive and the multiplicative notation for divisors and line bundles
will be mixed sometimes.
If $E$ is a vector bundle on $X$ and $L\in \mbox{Pic}(X)$, we set $E(-L):=E\otimes L^*$;
this notation will be used especially when $E$ is
replaced by the canonical bundle $K_C$ of a curve $C$.

\section{Definition, Properties, the First Applications}
\label{section: Def}

\subsection{Definition and First Properties}
\label{subsection: Definition LM}

We fix $S$ a smooth, complex, projective $K3$ surface and $L$ a globally
generated line bundle on $S$ with
$L^2=2g-2$. Let $C\in|L|$ be a
smooth curve and $A$ be a base-point-free line bundle
in $W^r_d(C)\setminus W^{r+1}_d(C)$.
As mentioned in the Introduction, the definition of Lazarsfeld--Mukai
bundles emerged from the attempt to lift
the linear system
$A$ to the surface $S$. Since it is virtually impossible to
lift it to another linear system, a higher-rank vector bundle
is constructed such that $H^0(C,A)$ corresponds to an $(r+1)$-dimensional
space of global sections. Hence $|A|$ lifts to a higher-rank
analogue of a linear system.

The kernel of the evaluation of sections of $A$\vspace*{-3pt}
\begin{equation}
\label{eqn: F}
 0\to F_{C,A}\to H^0(C,A)\otimes \mathcal{O}_S \buildrel{\mathrm{ev}}\over{\to} A\to 0\vspace*{-3pt}
\end{equation}
is a vector bundle of rank $(r+1)$.

\begin{defn}[Lazarsfeld \cite{LazarsfeldJDG}, Mukai \cite{MukaiPNAS}]
The \textit{Lazarsfeld--Mukai bundle} $E_{C,A}$ associated to the pair
$(C,A)$ is the dual of $F_{C,A}$.
\end{defn}

By dualizing the sequence~(\ref{eqn: F}) we obtain the short exact sequence\vspace*{-3pt}
\begin{equation}
\label{eqn: E}
 0\to H^0(C,A)^*\otimes \mathcal{O}_S \to E_{C,A}\to K_C(-A)\to 0,\vspace*{-3pt}
\end{equation}
and hence $E_{C,A}$ is obtained
from the trivial bundle by modifying it along the curve $C$
and comes equipped with a natural $(r+1)$-dimensional
space of global sections as planned.

We note here the first properties of $E_{C,A}$:

\begin{prop}[Lazarsfeld]
\label{prop: E_{C,A}}
The invariants of $E$ are the following:
\begin{itemize}
 \item[(1)] $\det (E_{C,A})= L$.
 \item[(2)] $c_2(E_{C,A})=d$.
 \item[(3)] $h^0 (S,E_{C,A})= h^0(C,A)+h^1(C,A)= 2r-d+1+g$.
 \item[(4)] $h^1(S,E_{C,A})=h^2(S,E_{C,A})=0$.
 \item[(5)] $\chi(S,E_{C,A}\otimes F_{C,A})=2(1-\rho(g,r,d))$.
 where $\rho(g,r,d)=g-(r+1)(g-d+r)$.
 \item[(6)] $E_{C,A}$ is globally generated off the base locus of $K_C(-A)$;
 in particular, $E_{C,A}$ is globally generated if $K_C(-A)$
 is globally generated.
\end{itemize}
\end{prop}
It is natural to ask conversely if given $E$ a vector
bundle on $S$ with  $\mbox{rk}(E)=r+1$, $h^1(S,E)=h^2(S,E)=0$, and
$\det (E)= L$, $E$ is the Lazarsfeld--Mukai bundle associated to a pair
$(C,A)$.
To this end, note that there is a rational map\vspace*{-3pt}
\[
 h_E : G(r+1,H^0(S,E)) \dashrightarrow |L|\vspace*{-3pt}
\]
defined in the following way. A general subspace $\varLambda\in G(r+1,H^0(S,E))$
is mapped to the degeneracy locus of the evaluation map:
$
\mathrm{ev} _{\varLambda} : \varLambda \otimes \mathcal{O}_S \to E.
$
If~the image $h_E(\varLambda)$ is a smooth curve $C\in |L|$, we set
$\mbox{Coker}(\mathrm{ev}_{\varLambda}):=K_C(-A)$,
where $A\in \mbox{Pic}(C)$ and $\mbox{deg}(A)=c_2(E)$, and
observe that $E=E_{C,A}$. Indeed, since $h^1(S,E)=0$, $A$
is globally generated, and from $h^2(S,E)=0$ it follows
that $\varLambda\cong H^0(C,A)^*$.
The conclusion is that:

\begin{prop}
\label{prop: caracterizare LM}
A rank-$(r+1)$ vector bundle $E$ on $S$ is a Lazarsfeld--Mukai bundle if
and only if $H^1(S,E)=H^2(S,E)=0$ and  there exists an $(r+1)$-dimensional subspace
of sections $\varLambda\subset H^0(S,E)$, such that the degeneracy locus of
the morphism $\mathrm{ev}_{\varLambda}$ is a smooth curve. In particular, being a
Lazarsfeld--Mukai vector bundle is an \textit{open condition}.
\end{prop}

Note that
there might be different pairs with the same Lazarsfeld--Mukai
bundles, the difference being given by the corresponding spaces
of global sections.

\subsection{Simple and Non-simple Lazarsfeld--Mukai Bundles}
\label{subsection: Simple LM}

We keep the notation from the previous subsection.
In the original situation, the bundles used by Lazarsfeld \cite{LazarsfeldJDG}
and Mukai \cite{MukaiPNAS} are simple.
The non-simple Lazarsfeld--Mukai bundles are, however, equally useful
\cite{Aprodu-FarkasCOMP,Ciliberto-PareschiCRELLE}.
For instance, Lazarsfeld's argument is partly based on an analysis of the
non-simple bundles.

Proposition~\ref{prop: E_{C,A}} already shows that for $\rho(g,r,d)<0$
the associated Lazarsfeld--Mukai bundle cannot be simple.
The necessity of making a distinction between simple and non-simple bundles
for nonnegative $\rho$ will become more evident in the next sections.

In the rank-two case, one can give a precise description \cite{Donagi-MorrisonJDG} of non-simple
Lazarsfeld--Mukai bundles, see also \cite{Ciliberto-PareschiCRELLE} Lemma 2.1:

\begin{lem}[Donagi--Morrison]
\label{lemma: DM}
Let $E_{C,A}$ be a non-simple Lazarsfeld--Mukai bundle.
Then there exist line bundles $M,N\in \mathrm{Pic}(S)$ such that
$h^0(S,M)$, $h^0(S,N)\ge 2$,
$N$ is globally generated, and there exists a
locally complete intersection subscheme $\xi$ of $S$,
either of dimension zero or the empty set, such that
$E_{C,A}$ is expressed as an extension
\begin{equation}
\label{eq: DM}
0\to M\to E_{C,A}\to N\otimes I_{\xi} \to 0.
\end{equation}
Moreover, if $h^0(S,M\otimes N^*)=0$, then $\xi=\emptyset$ and the extension splits.
\end{lem}

One can prove furthermore that $h^1(S,N)=0$, \cite{Aprodu-FarkasCOMP} Remark 3.6.

We say that~(\ref{eq: DM}) is the \emph{Donagi--Morrison extension}
associated to $E_{C, A}$. This notion makes perfect sense as
this extension is uniquely determined by the vector bundle,
if it is indecomposable \cite{Aprodu-FarkasCOMP}.
Actually, a \textit{decomposable} Lazarsfeld--Mukai bundle $E$ cannot
be expressed as an extension~(\ref{eq: DM}) with $\xi\ne\emptyset$, and
hence a Donagi--Morrison extension is always
unique, up to a permutation of factors in the decomposable case.
Moreover, a Lazarsfeld--Mukai bundle is decomposable
if and only if the corresponding Donagi--Morrison extension is trivial.

In the higher-rank case, we do not have such a precise description.\footnote{In fact,
we do have a Harder--Narasimhan  filtration,
but we cannot control all the factors.} However, a similar
sufficiently strong statement is still valid \cite{LazarsfeldJDG,LazarsfeldICTP,PareschiJAG}.

\begin{prop}[Lazarsfeld]
 Notation as above. If $E_{C,A}$ is not simple, then 
the linear system $|L|$ contains a reducible or a multiple
curve.
\end{prop}

In the rank-two case, this statement comes from the decomposition
$L\cong M\otimes N$.

\subsection{The Petri Conjecture Without Degenerations}
\label{subsection: BNP}

A smooth curve of genus $g$ is said to satisfy \textit{Petri's condition},
or to be \textit{Brill--Noether--Petri generic},
if the multiplication map (the Petri map)
\[
 \mu_{0,A}: H^0(C, A) \otimes H^0(C, K_C(-A))
 \to  H^0(C, K_C),
\]
is injective for any line bundle $A$ on $C$.
One consequence of this condition is that
all the Brill--Noether loci $W^r_d(C)$ have the expected
dimension and are smooth away from $W^{r+1}_d(C)$;
recall that the tangent space at the point $[A]$
to $W^r_d(C)$ is naturally isomorphic to the dual
of $\mathrm{Coker}(\mu_{0,A})$ \cite{Arbarello-Cornalba-Griffiths-Harris}.
The Petri conjecture, proved by degenerations by Gieseker,
states that a general curve satisfies Petri's condition.
Lazarsfeld \cite{LazarsfeldJDG} found a simpler and
elegant proof without degenerations by analyzing curves
on very general $K3$ surfaces.

Lazarsfeld's idea is to relate the Petri maps
to the Lazarsfeld--Mukai bundles; this relation is valid
in general and has many other
applications. Suppose, as in the previous subsections,
that $S$ is a $K3$ surface and $L$ is a globally generated
line bundle on $S$. For the moment, we do not need to assume
that $L$ generates the Picard group.
E. Arbarello and M. Cornalba  constructed
a scheme $\mathcal{W}^r_d(|L|)$
parameterizing pairs $(C,A)$ with $C\in |L|$ smooth and
$A\in W^r_d(C)$ and a morphism\vspace*{-3pt}
\[
\pi_S:\mathcal{W}^r_d(|L|)\to |L|.\vspace*{-3pt}
\]

Assume that $A\in W^r_d(C)\setminus W^{r+1}_d(C)$ is globally generated,
and consider $M_A$ the vector bundle of rank $r$
on $C$ defined as the kernel of the evaluation map
\begin{equation}\label{MA}
 0\to M_A \to  H^0(C,A)\otimes \mathcal{O}_C \buildrel{\mathrm{ev}}\over{\to} A\to 0.
\end{equation}
Twisting~(\ref{MA}) with $K_C\otimes A^*$, we obtain the
following description of the kernel of the Petri map:\footnote{This ingenious procedure is an efficient replacement
of the base-point-free pencil trick;
``it has killed the base-point-free pencil trick,'' to quote Enrico Arbarello.}\vspace*{4pt}
\[
 \mbox{Ker}(\mu_{0,A})= H^0(C, M_A\otimes K_C\otimes A^*).
\]

There is another exact sequence on $C$\vspace*{-4pt}
\[
0\to \mathcal{O}_C\to F_{C,A}|_C\otimes K_C\otimes A^*\to M_A\otimes K_C\otimes A^*\to 0,\vspace*{-4pt}
\]
and from the  defining sequence of $E_{C,A}$ one obtains the exact sequence on $S$\vspace*{-4pt}
\[
0\to H^0(C, A)^*\otimes F_{C,A}\to E_{C,A}\otimes F_{C,A}\to
F_{C, A}|_C\otimes K_C\otimes A^*\to 0.\vspace*{-4pt}
\]

From the vanishing of $h^0(C,F_{C,A})$
and of $h^1(C, F_{C, A})$, we obtain\vspace*{-4pt}
\[
H^0(C,E_{C,A}\otimes F_{C,A})=H^0(C,F_{C,A}|_C\otimes K_C\otimes A^*).\vspace*{-4pt}
\]

Suppose that $\mathcal{W}\subset \mathcal{W}^r_d(|L|)$ is a dominating component and
$(C,A)\in\mathcal{W}$ is an element such that $A$ is
globally generated and $h^0(C,A)=r+1$. A deformation-theoretic
argument shows that if the Lazarsfeld--Mukai
bundle $E_{C,A}$ is simple, then the coboundary map
$H^0(C,M_A\otimes K_C\otimes A^*)\to H^1(C,\mathcal{O}_C)$
is zero \cite{PareschiJAG}, which eventually implies the
injectivity of $\mu_{0,A}$.

By reduction to complete base-point-free bundles on the curve \cite{LazarsfeldJDG,PareschiJAG} this analysis yields:

\begin{thm}[Lazarsfeld]
\label{thm: Lazarsfeld}
Let $C$ be a smooth curve of genus $g\ge 2$ on a $K3$
surface~$S$, and assume that any divisor in the linear
system $|C|$ is reduced and irreducible. Then a
generic element in the linear system $|C|$
is Brill--Noether--Petri generic.
\end{thm}

A particularly interesting case is when the Picard group of $S$ is
generated by $L$ and $\rho(g,r,d)=0$. Obviously, the condition
$\rho=0$ can be realized only for composite genera, as
$g=(r+1)(g-d+r)$, for example, $r=1$ and $g$ even.
Under these assumptions, there is a unique
Lazarsfeld--Mukai bundle $E$ with $c_1(E)=L$ and $c_2(E)=d$,
and different pairs $(C,A)$ correspond to different
\hbox{$\varLambda\in G(r+1,$} $H^0(S,E))$; in other words
the natural rational map $G(r+1,H^0(S,E))\dashrightarrow
\mathcal{W}^r_d(|L|)$ is dominating.
Note that $E$ must be stable and globally generated.

\subsection{Mukai Manifolds of Picard Number One}
\label{subsection: Mukai manifolds}

A Fano manifold $X$ of dimension $n\ge 3$
and index $n-2$ (i.e., of coindex 3) is called
a \textit{Mukai manifold}.\footnote{Some authors consider
that Mukai manifolds have dimension four or more.}
In the classification, special attention is given to prime Fano manifolds:
note that if $n\ge 7$, $X$ is automatically prime as shown by Wisniewski;
see, for example, \cite{Iskovskih-Prokhorov}.

Assume  that the Picard group of $X$ is generated
by an ample line bundle $L$, and let the sectional genus
$g$ be the integer $(L^n)/2+1$. Mukai and Gushel
used vector bundle techniques to obtain a complete
classification of these manifolds. A~first major obstacle
is to prove that the fundamental linear system contains indeed
a smooth element, aspect which is settled by
Shokurov and Mella; see, for example, \cite{Iskovskih-Prokhorov}.
Then the $(g+n-2)$-dimensional linear system
$|L|$ is base-point-free, and
a general linear section
with respect to the generator of the Picard group
is a $K3$ surface. More precisely, if $\mathrm{Pic}(X)
=\mathbb{Z}\cdot L$, then for $H_1,\cdots,H_{n-2}$
general elements in the fundamental linear system
$|L|$, $S:=H_1\cap\cdots\cap H_{n-2}$ is
scheme-theoretically a $K3$ surface. Note that if $n\ge 4$
and $i\ge 3$, the intersection $H_1\cap\cdots\cap H_{n-i}$
is again a Fano manifold of coindex 3.

Mukai noticed that the fundamental linear system either
is very ample, and the image of $X$ is projectively
normal or is associated to a double covering of $\mathbb{P}^n$
($g=2$) or of the hyper-quadric $Q^n\subset\mathbb{P}^{n+1}$
($g=3$). The difficulty of the problem is thus
to classify all the possible cases where
$|L|$ is normally generated, called \textit{of the first
species}. Taking linear
sections one reduces (not quite immediately) to the case $n=3$
\cite{Iskovskih-Prokhorov} p.110.

For simplicity, let us assume that
$X$ is a prime Fano 3-fold of index 1. If
$g=4$ and $g=5$, $X$ is a complete intersection; hence
the hard cases begin with genus $6$.
A hyperplane section $S$ is a $K3$ surface, and,
by a result of Moishezon, $\mathrm{Pic}(S)$ is generated
by $L|_S$.

Let us denote by $\mathcal{F}_g$ the moduli
space of polarized $K3$ surfaces of degree $2g-2$, by $\mathcal{P}_g$
the moduli space of pairs $(K3\ \mbox{surface},\mbox{curve})$ and
$\mathcal{M}_g$ the moduli space of genus-$g$ curves.
There are two nice facts in Mukai's proof involving
these two moduli spaces. His first observation is
that if there exists a prime Fano 3-fold $X$
of the first species of genus $g\ge 6$ and index $1$,
the rational map $\phi_g:\mathcal{P}_g\dashrightarrow \mathcal{M}_g$
is \textit{not} generically finite \cite{MukaiLMS}. The second nice
fact is that $\phi_g$ is generically finite
if and only if $g=11$ or $g\ge 13$ \cite{MukaiLMS}.\footnote{In genus $11$, it
is actually birational \cite{MukaiLNPAM}.} Hence, one
is reduced to study the genera $6\le g\le 12$
with $g\ne 11$. At this point, Lazarsfeld--Mukai
bundles are employed. By the discussion from Sect.\,\ref{subsection: BNP},
for any decomposition $g=(r+1)(g-d+r)$, with
$r\ge 1$, $d\le g-1$, there exists a unique
Lazarsfeld--Mukai bundle $E$ of rank $(r+1)$.
It has already been noticed that the bundle $E$ is stable and globally generated.
Moreover, the determinant map\vspace*{-4pt}
\[
 \mathrm{det}:\wedge^{r+1}H^0(S,E)\to H^0(S,L)\vspace*{-4pt}
\]
is surjective \cite{MukaiPNAS}, and hence it induces
a linear embedding\vspace*{-4pt}
\[
 \mathbb{P}H^0(S,L)^*\hookrightarrow \mathbb{P}(\wedge^{r+1}H^0(S,E)^*).\vspace*{-4pt}
\]

Following \cite{MukaiPNAS}, we have a commutative diagram
\begin{center}
\mbox{\xymatrix{S \ar[r]^{\phi_E}\ar@{^(->}[d]^{\phi_{|L|}} & G \ar@{^(->}[d]^{\mathrm{Pluecker}}\\
\mathbb{P}H^0(L)^*\ar@{^(->}[r] & \mathbb{P}(\wedge^{r+1}H^0(E)^*)}}
\end{center}
where $G:=G(r+1,H^0(S,E)^*)$ and $\phi_E$ is given by $E$.
This diagram shows that $S$ is embedded in a suitable linear section
of the Grassmannian $G$. Moreover, this diagram
extends over $X$: by a result of Fujita, $E$
extends to a stable vector bundle on $X$, and
the diagram over $X$ is obtained for similar reasons.
Hence $X$ is a linear section of a Grassmannian.
By induction on the dimension, $X$ is contained in
a \textit{maximal} Mukai manifold, which is also
a linear section of the Grassmannian.
A complete list of maximal Mukai manifolds
is given in \cite{MukaiPNAS}. Notice that
in genus $12$, the maximal Mukai manifolds are
threefold already.

\setcounter{thm}{0}
\section{Constancy of Invariants of $K3$ Sections}
\label{section: Constancy}

\subsection{Constancy of the Gonality. I}
\label{subsection Gonality I}

In his analysis of linear systems on $K3$ surfaces Saint--Donat
\cite{Saint-DonatAJM}
shows that any smooth curve which is linearly equivalent
to a hyperelliptic or trigonal curve is also hyperelliptic, respectively trigonal.
The idea was to prove that the
minimal pencils are induced by elliptic pencils
defined on the surface. This result was sensibly extended by Reid
 \cite{ReidJLM} who proved the following existence result:

\begin{thm}[Reid]
Let $C$ be a smooth curve of genus $g$ on a $K3$ surface $S$
and $A$ be a complete, base-point-free $g^1_d$ on
$C$. If\vspace*{-3pt}
\[
\frac{d^2}{4}+d+2<g,\vspace*{-3pt}
\]
then $A$ is the restriction of an elliptic pencil on $S$.
\end{thm}

It is a good occasion to present here, as a direct application of techniques involving
Lazarsfeld--Mukai bundles, an alternate shorter proof of Reid's theorem.

\begin{proof}
We use the notation of previous sections.
By the hypothesis, the Lazarsfeld--Mukai bundle $E$ is not simple, and
hence we have a unique Donagi--Morrison extension\vspace*{3pt}
\[
0\to M\to E\to N\otimes I_{\xi}\to 0,\vspace*{2pt}
\]
with $\xi$ of length $\ell$. Note that $M\cdot N=d-\ell\le d$.
By the Hodge index theorem, we have
$(M^2)\cdot(N^2)\le (M\cdot N)^2\le d^2$, whereas from
$M+N=C$ we obtain $(M^2)=2(g-1-d)-(N^2)$, hence\vspace*{-4pt}
\[
(N^2)\le\frac{d^2}{2(g-1-d)-(N^2)}.\vspace*{-4pt}
\]

Therefore, the even integer $x:=(N^2)$ satisfies the following
inequality
$x^2-2x(g-1-d)+d^2\ge 0.$
The hypothesis shows that the above inequality fails for $x\ge 2$,
and hence $N$ must be an elliptic pencil.
\end{proof}

In conclusion, for small values, the gonality\footnote{The gonality
$\mathrm{gon}(C)$ of a curve $C$ is the minimal degree of a morphism
from $C$ to the projective~line.}
is constant in the linear system. Motivated
by these facts, Harris and Mumford conjectured
that \textit{the gonality of \hbox{$K3$-sections}
should always be constant} \cite{Harris-MumfordINVENTIONES}.

This conjecture is unfortunately wrong as stated: Donagi and Morrison
\cite{Donagi-MorrisonJDG} gave the following counterexample:

\begin{ex}
\label{ex: DM}
 Let $S\to \mathbb{P}^2$ be a double cover branched
 along a smooth sextic and $L$ be the pull-back of
 $\mathcal{O}_{\mathbb{P}^2}(3)$. The curves in $|L|$
 have all genus $10$. The general curve $C\in|L|$
 is isomorphic to a smooth plane sextic, and hence
 it is pentagonal. On the other hand, the pull-back of
 a general smooth plane cubic $\varGamma$ is a double
 cover of $\varGamma$, and thus it is tetragonal.
\end{ex}

\subsection{Constancy of the Clifford Index}
\label{subsection: Cliff}

Building on his work on Koszul cohomology and its relations
with geometry, M. Green proposed a reformulation of the
Harris-Mumford conjecture replacing the gonality by
the Clifford index.

\medskip

Recall that the \textit{Clifford index} of
a nonempty linear system $|A|$ on a smooth curve $C$
is the codimension of the image
of the natural addition map $|A|\times|K_C(-A)|\to |K_C|$.
This definition is nontrivial only for relevant linear
systems $|A|$, i.e., such that both $|A|$ and $|K_C(-A)|$
are at least one-dimensional; such an $A$ is said
to \textit{contribute to the Clifford index}.
The \textit{Clifford index of $C$} is the minimum of all
the Clifford indices taken over the  linear
systems that contribute to the Clifford index and
is denoted by $\mathrm{Cliff}(C)$.
The Clifford index is related to the gonality
by the following inequalities\vspace*{3pt}
\[
\mathrm{gon}(C)-3\le \mathrm{Cliff}(C)\le \mathrm{gon}(C)-2,\vspace*{3pt}
\]
and curves with $\mathrm{gon}(C)-3 = \mathrm{Cliff}(C)$
are very rare: typical examples are plane curves and
Eisenbud--Lange--Martens--Schreyer curves \cite{ELMS,KnutsenIJM}.\footnote{It
is conjectured that the only other examples should
be some half-canonical curves of even genus
and maximal gonality \cite{ELMS}; however,
this conjecture seems to be very difficult.}

From the Brill--Noether theory, we obtain the
bound $\mathrm{Cliff}(C)\le\left[(g-1)/2\right]$
(and, likewise, $\mathrm{gon}(C)\le\left[(g+3)/2\right]$),
and it is known that the equality is achieved for general curves.
The Clifford index is in fact a measure of how special
a curve is in the moduli space.

\medskip

The precise statement obtained by Green and Lazarsfeld
is the following \cite{Green-LazarsfeldINVENTIONES}:

\begin{thm}[Green--Lazarsfeld]
\label{thm: GL Cliff}
 Let $S$ be a $K3$ surface and $C\subset S$ be a smooth
 irreducible curve of genus $g\ge 2$. Then $\mathrm{Cliff}(C')
 =\mathrm{Cliff}(C)$ for every smooth curve $C'\in|C|$. Furthermore,
 if $\mathrm{Cliff}(C)$ is strictly less than the generic value
 $ \left[(g-1)/2\right]$, then there exists a line bundle
 $M$ on $S$ whose restriction to any smooth curve $C'\in|C|$ computes
 the Clifford index of~$C'$.
\end{thm}

The proof strategy is  based on a reduction method of the
associated Lazarsfeld--Mukai bundles. The bundle $M$
is obtained from the properties of the reductions; we
refer to \cite{Green-LazarsfeldINVENTIONES} for
details.

From the Clifford index viewpoint, Donagi--Morrison's
example is not different from the other cases. Indeed, all
smooth curves in $|L|$ have Clifford index $2$.
We shall see in the next subsection that
Donagi--Morrison's example
is truly an isolated exception for the constancy of the gonality.

\subsection{Constancy of the Gonality. II}
\label{subsection: Gonality II}

As discussed above, the Green--Lazarsfeld proof of the
constancy of the Clifford index was mainly based on
the analysis of Lazarsfeld--Mukai bundles.
It is natural to try and explain the peculiarity of
Donagi--Morrison's example from this point of view.
This was done in \cite{Ciliberto-PareschiCRELLE}.
The surprising answer found by Ciliberto and
Pareschi \cite{Ciliberto-PareschiCRELLE}
(see also \cite{Donagi-MorrisonJDG})
is the following:

\begin{thm}[Ciliberto--Pareschi]
\label{thm: Ciliberto-Pareschi}
 Let $S$ be a $K3$ surface and $L$ be an ample line bundle on $S$.
If the gonality of the smooth curves in $|L|$ is not constant,
then $S$ and $L$ are as in Donagi--Morrison's example.
\end{thm}

Theorem~\ref{thm: Ciliberto-Pareschi} was refined
by Knutsen \cite{KnutsenIJM} who replaced ampleness by the more
general condition that $L$ be globally generated.
The extended setup covers also the case of exceptional curves, as introduced
by Eisenbud, Lange, Martens, and Schreyer \cite{ELMS}.

The proof of Theorem~\ref{thm: Ciliberto-Pareschi}
consists of a thorough analysis of the
loci $\mathcal{W}^1_d(|L|)$, where $d$ is the minimal
gonality of smooth curves in $|L|$, through the
associated Lazarsfeld--Mukai bundles. The authors identify
Donagi--Morrison's example in the following way:

\begin{thm}[Ciliberto--Pareschi]
\label{thm: Ciliberto-Pareschi2}
Let $S$ be a $K3$ surface and $L$ be an ample line bundle on $S$.
If the gonality of smooth curves in $|L|$ is not constant and if
there is a pair $(C,A)\in\mathcal{W}^1_d(|L|)$ such that
$h^1(S,E_{C,A}\otimes F_{C,A})=0$, then $S$ and
$L$ are as in Donagi--Morrison's example.
\end{thm}

To conclude the proof of Theorem~\ref{thm: Ciliberto-Pareschi}, Ciliberto and Pareschi prove that
non-constancy of the gonality implies the existence
of a pair $(C,A)$ with $h^1(S,E_{C,A}\otimes F_{C,A})=0$;
see \cite{Ciliberto-PareschiCRELLE} Proposition 2.4.

It is worth to notice that, in Example~\ref{ex: DM}, if $C$ is the
inverse image of a plane cubic and $A$ is a $g^1_4$
(the pull-back of an involution), then
$E_{C,A}$ is the pull-back of
$\mathcal{O}_{\mathbb{P}^2}(1)\oplus\mathcal{O}_{\mathbb{P}^2}(2)$
\cite{Ciliberto-PareschiCRELLE}, and hence
the vanishing of  $h^1(S,E_{C,A}\otimes F_{C,A})$ is guaranteed
in this case.

\subsection{Parameter Spaces of Lazarsfeld--Mukai
Bundles and~Dimension of  Brill--Noether Loci}
\label{subsection: Parameter}

We have already seen that the Brill--Noether loci
are smooth of expected dimension at pairs corresponding
to simple Lazarsfeld--Mukai bundles. It
is interesting to know what is the dimension of
these loci at other points as well. Precisely,
we look for a uniform bound on the dimension
of Brill--Noether loci of general curves in a linear
system.

A first step was made by Ciliberto and Pareschi
\cite{Ciliberto-PareschiCRELLE} who proved,
as a necessary step in Theorem~\ref{thm: Ciliberto-Pareschi},
that an ample curve of gonality strictly less
than the generic value, general in its linear system,
carries finitely many minimal pencils. This result
was  extended to other Brill--Noether loci
\cite{Aprodu-FarkasCOMP},
proving a phenomenon of \textit{linear growth} with the degree; see below.
Let us mention that, for the moment, the only results
in this direction are known to hold for pencils
\cite{Aprodu-FarkasCOMP} and nets \cite{Lelli-Chiesa}.

As before, we consider $S$ a $K3$ surface and $L$
a globally generated line bundle on $S$.
In order to parameterize all pairs $(C, A)$
with non-simple  Lazarsfeld--Mukai bundles, we need a global construction.
We fix a nontrivial globally generated line bundle
$N$ on $S$ with $H^0(L(-2N))\neq 0$
and an integer $\ell\ge 0$. We set $M:=L(-N)$
and $g:=1+L^2/2$.
Define $\widetilde{\mathcal{P}}_{N,\ell}$ to be the family of \textit{vector bundles} of rank $2$ on
 $S$ given by nontrivial extensions\vspace*{-3pt}
\begin{equation}
\label{eq: extension}
0\to M\to E\to N\otimes I_{\xi}\to 0,\vspace*{-3pt}
\end{equation}
where $\xi$ is a zero-dimensional locally complete intersection subscheme
(or the empty set) of $S$ of length $\ell$, and set
\[
\mathcal{P}_{N,\ell}:=\{[E]\in\widetilde{\mathcal{P}}_{N, \ell}:\ h^1(S,E)=h^2(S,E)=0\}.
\]
Equivalently (by Riemann--Roch), $[E]\in \mathcal{P}_{N, \ell}$ if and only if
$h^0(S,E)=$\break $g-c_2(E)+3$ and $h^1(S,E)=0$. Note that any non-simple
Lazarsfeld--Mukai bundle on $S$ with
determinant $L$ belongs to some family $\mathcal{P}_{N,\ell}$,
from Lemma~\ref{lemma: DM}.
The family $\mathcal{P}_{N,\ell}$,
which, a priori, might be the empty set, is an open Zariski subset
of a projective bundle of the Hilbert scheme~$S^{[\ell]}$.

Assuming that $\mathcal{P}_{N,\ell}\ne \emptyset$,
we consider  the  Grassmann bundle $\mathcal{G}_{N,\ell}$ over $\mathcal{P}_{N,\ell}$
classifying  pairs $(E,\varLambda)$ with $[E]\in\mathcal{P}_{N,\ell}$ and
$\varLambda\in \mathrm{G}(2,H^0(S,E))$. If $d:=c_2(E)$ we define the rational map $h_{N, \ell}: \mathcal{G}_{N, \ell}
\dashrightarrow \mathcal{W}^1_d(|L|)$, by setting $h_{N, \ell}(E, \varLambda):=(C_{\varLambda}, A_{\varLambda})$, where $A_{\varLambda}\in \mbox{Pic}^d(C_{\varLambda})$ is such that the following exact sequence on $S$ holds:\vspace*{4pt}
\[
0\to \varLambda\otimes \mathcal{O}_S\stackrel{\mathrm{ev}_{\varLambda}}\to E\to K_{C_{\varLambda}}\otimes A_{\varLambda}^*\to 0.\vspace*{3pt}
\]

One computes
$\dim\ \mathcal{G}_{N,\ell}=g+\ell+h^0(S,M\otimes N^*)$.
If we assume furthermore that $\mathcal{P}_{N,\ell}$ contains a Lazarsfeld--Mukai vector bundle $E$ on $S$ with $c_2(E)=d$
and consider $\mathcal{W}\subset \mathcal{W}^1_d(|L|)$ the closure of the
image of the rational map $h_{N,\ell}:\mathcal{G}_{N,\ell}\dashrightarrow \mathcal{W}^1_d(|L|)$,
then we find $\dim\ \mathcal{W} =g+d-M\cdot N=g+\ell$.

On the other hand, if $C\in|L|$ has Clifford
dimension one and $A$ is a globally generated line
bundle on $C$ with $h^0(C,A)=2$ and $[E_{C,A}]\in\mathcal{P}_{N,\ell}$,
then $M\cdot N\ge \mathrm{gon}(C)$.

These considerations on the indecomposable
case, together with a simpler analysis of
decomposable bundles, yield finally \cite{Aprodu-FarkasCOMP}:

\begin{thm}
\label{thm: Green Cliffdim 1}
Let $S$ be a $K3$ surface and $L$ a globally generated line bundle on $S$,
such that general curves in $|L|$ are of Clifford dimension one. Suppose that
$\rho(g,1,k)\leq 0$, where $L^2=2g-2$ and $k$ is the (constant) gonality of
all smooth curves in $|L|$. Then for a
general curve $C\in|L|$, we have
\begin{equation}
\label{lgc}
  \mathrm{dim}\ W^1_{k+d}(C)=d\mbox{ for all }
0\le d\le g-2k+2.
\end{equation}
\end{thm}

The condition~(\ref{lgc}) is called the \textit{linear growth condition}.
It is equivalent to
\[
\mathrm{dim}\ W^1_{g-k+2}(C)=\rho(g,1,g-k+2)=g-2k+2.
\]

Note that the condition that $C$ carry finitely many minimal
pencils, which is a part of~(\ref{lgc}), appears explicitly in~\cite{Ciliberto-PareschiCRELLE}.
It is directly related to the constancy of the gonality
discussed before.

\setcounter{thm}{0}
\section{Green's Conjecture for Curves on $K3$ Surfaces}
\label{section: Green}

\subsection{Koszul Cohomology}\label{ch1:sec3.1}

Let $X$ be a (not necessarily smooth) complex, irreducible,
projective variety and $L\in\mathrm{Pic}(X)$ globally generated.
The Euler sequence on the projective space $\mathbb{P}(H^0(X,L)^*)$
pulls back to a short exact sequence of vector bundles
on~$X$
\begin{equation}
\label{eqn: Euler}
0 \to M_L\to H^0(X,L)\otimes \mathcal{O}_X\to L\to 0.
\end{equation}

After taking exterior powers in the sequence~(\ref{eqn: Euler}),
twisting with multiples of $L$ and going to global sections,
we obtain an exact sequence for any nonnegative $p$ and $q$:\vspace*{-5pt}
\begin{align}
\label{eqn: WedgeEuler}
0\to H^0(\wedge^{p+1}M_L\otimes L^{q-1})
\to \wedge^{p+1}H^0(L)\otimes H^0(L^{q-1})
\stackrel{\delta}{\to} H^0(\wedge^pM_L\otimes L^q).\nonumber\\ \vspace*{-5pt}
\end{align}

The finite-dimensional vector space $K_{p,q}(X,L):=\mathrm{Coker}(\delta)$
is called the \textit{Koszul cohomology space}\footnote{The
indices $p$ and $q$ are usually forgotten when defining Koszul cohomology.}
of $X$ with values in $L$
\cite{LazarsfeldICTP,GreenJDG,GreenICTP}.
Observe that $K_{p,q}$ can be defined alternatively as:
\[
K_{p,q}(X,L)=\mathrm{Ker}\left(H^1(\wedge^{p+1}M_L\otimes L^{q-1})
\to \wedge^{p+1}H^0(L)\otimes H^1(L^{q-1})\right),\vspace*{-3pt}
\]
description which is particularly useful when $X$ is a curve.

Several versions are used in practice, for example, replace
$H^0(L)$ in~(\ref{eqn: Euler}) by a subspace that generates $L$
or twist~(\ref{eqn: WedgeEuler}) by $\mathcal{F}\otimes L^{q-1}$
where $\mathcal{F}$ is a coherent sheaf. For our presentation,
however, we do not need to discuss these natural generalizations.

Composing the maps\vspace*{-3pt}
\[
\wedge^{p+1}H^0(L)\otimes H^0(L^{q-1})
\stackrel{\delta}{\to} H^0(\wedge^pM_L\otimes L^q)
\hookrightarrow \wedge^pH^0(L)\otimes H^0(L^q)\vspace*{-3pt}
\]
we obtain, by iteration,  a complex\vspace*{-3pt}
\[
\wedge^{p+1}H^0(L)\otimes H^0(L^{q-1})\to
\wedge^pH^0(L)\otimes H^0(L^q)\to
\wedge^{p-1}H^0(L)\otimes H^0(L^{q+1})\vspace*{-3pt}
\]
whose cohomology at the middle  is
$K_{p,q}(X,L)$, and this is the definition given by Green
\cite{GreenJDG}.

\medskip

An important property of Koszul cohomology is upper-semicontinuity
in flat families with constant cohomology; in particular,
vanishing of Koszul cohomology is an open
property in such families. For curves, constancy
of $h^1$ is a consequence of flatness and of constancy of $h^0$, as
shown by the Riemann--Roch theorem.

\medskip

The original motivation for studying Koszul cohomology
spaces was given by the relation with minimal resolutions
over the polynomial ring. More precisely, if $L$ is very
ample, then the Koszul cohomology computes the minimal resolution
of the graded module\vspace*{4pt}
\[
R(X,L):=\bigoplus_qH^0(X,L^q)\vspace*{3pt}
\]
over the polynomial ring
\cite{GreenJDG,GreenICTP}; see also \cite{EisenbudBOOK,Aprodu-NagelULECT}, in the sense that any graded piece
that appears in the minimal resolution is (non-canonically) isomorphic to a
$K_{p,q}$. If the image of $X$ is projectively normal,
this module coincides with the homogeneous coordinate
ring of $X$. The projective normality of $X$ can also
be read off Koszul cohomology, being
characterized by the vanishing condition $K_{0,q}(X,L)=0$ for all
$q\ge 2$. Furthermore, for a projectively normal $X$,
the homogeneous ideal is generated by quadrics if and
only if $K_{1,q}(X,L)=0$ for all $q\ge 2$.\footnote{The
dimension of $K_{1,q}$ indicates the number of generators
of degree $(q+1)$ in the homogeneous~ideal.} The phenomenon
continues as follows: if $X$ is projectively normal and the homogeneous ideal
is generated by quadrics,
then the relations between the generators are linear if and only if
$K_{2,q}(X,L)=0$ for all $q\ge 2$, whence the relation with
syzygies ~\cite{GreenJDG}.

\medskip

Other notable application of Koszul cohomology is the description
of Castelnuovo--Mumford regularity, which coincides
with, \cite{GreenJDG,Aprodu-NagelULECT}
\[
\mathop{\mathrm{min}}_q\{K_{p,q}(X,L)=0,\mbox{ for all }p\}.
\]

Perhaps the most striking property of Koszul cohomology,
discovered by Green and Lazarsfeld \cite[Appendix]{GreenJDG},
is a consequence of a nonvanishing result:

\begin{thm}[Green--Lazarsfeld]
\label{thm: GL nonvan}
 Suppose $X$ is smooth and $L=L_1\otimes L_2$ with $r_i:=h^0(X,L_i)-1\ge 1$.
 Then $K_{r_1+r_2-1,1}(X,L)\ne 0$.
\end{thm}

Note that the spaces $K_{p,1}$ have the following particular
attribute: if $K_{p,1}\ne 0$ for some $p\ge 1$ then
$K_{p',1}\ne 0$ for all $1\le p'\le p$. This is obviously
false for $K_{p,q}$ with $q\ge 2$.

Theorem~\ref{thm: GL nonvan} shows that the existence of nontrivial decompositions
of $L$ reflects onto the existence of nontrivial
Koszul classes in some space $K_{p,1}$.
Its most important applications are for curves, in
particular for canonical curves, case which is discussed in the next subsection.
In the higher-dimensional cases, for surfaces, for instance, the
meaning of  Theorem~\ref{thm: GL nonvan} becomes more
transparent if it is accompanied by a restriction theorem which
compares the Koszul cohomology of $X$
with the Koszul cohomology of the linear sections
\cite{GreenJDG}:

\begin{thm}[Green]
\label{thm: Lefschetz}
 Suppose $X$ is smooth and $h^1(X,L^q)=0$ for all $q\ge 1$.
 Then for any connected reduced divisor $Y\in|L|$,
 the restriction map induces an isomorphism\vspace*{4pt}
 \[
 K_{p,q}(X,L)\stackrel{\sim}{\to}K_{p,q}(Y,L|_Y),\vspace*{3pt}
 \]
 for all $p$ and $q$.
\end{thm}

The vanishing of $h^1(X,\mathcal{O}_X)$ suffices
to prove that the restriction is an isomorphism between
the spaces $K_{p,1}$ \cite{Aprodu-NagelULECT}.

In the next subsections, we shall apply Theorem~\ref{thm: Lefschetz} for $K3$ sections.

\begin{cor}
\label{cor: Lefschetz K3}
Let $C$ be a smooth connected curve on a $K3$
surface $S$. Then\vspace*{-4pt}
\[
K_{p,q}(S,\mathcal{O}_S(C))\cong
K_{p,q}(C,K_C)\vspace*{-4pt}
\]
for all $p$ and $q$.
\end{cor}

One direct consequence is a duality theorem for 
Koszul cohomology of $K3$ surfaces.\footnote{Duality
for Koszul cohomology of curves follows from Serre's duality. For higher-dimensional
manifolds, some supplementary vanishing conditions are required \cite{GreenJDG,GreenICTP}.} It shows the
symmetry of the table containing the dimensions
of the spaces $K_{p,q}$, called \textit{the Betti
table}.

\vspace*{-4pt}
\subsection{Statement of Green's Conjecture}\label{ch1:sec3.2}

Let us particularize Theorem~\ref{thm: GL nonvan} for a
canonical curve. Consider $C$ a smooth curve and
choose a decomposition $K_C=A\otimes K_C(-A)$.
Theorem~\ref{thm: GL nonvan} applies only if
$h^0(C,A)\ge 2$ and $h^1(C,A)\ge 2$, i.e., if
$A$ contributes to the Clifford index.
The~quantity $r_1+r_2-1$ which appears in the statement
equals $g-\mathrm{Cliff}(A)-2$, and hence, if $A$
\textit{computes} the Clifford index, we obtain
the following:

\begin{thm}[Green--Lazarsfeld]
 For any smooth curve $C$ of genus $g$ Clifford index $c$
 we have $K_{g-c-2,1}(C,K_C)\ne 0$.
\end{thm}

It is natural to determine whether or not this result is sharp, question
which is addressed in the statement Green's conjecture:

\begin{con}[Green]
 Let $C$ be a smooth curve. For all $p\ge g-c-1$,
 we have $K_{p,1}(C,K_C) = 0$.
\end{con}

For the moment, Green's conjecture remains a hard open problem.
At the same time, strong evidence has been discovered.
For instance, it is known to hold for general curves \cite{VoisinJEMS,VoisinCOMP}, for curves of odd genus and maximal
Clifford index \cite{VoisinCOMP,Hirschowitz-RamananAENS},
for general curves of given gonality
\cite{VoisinJEMS,TeixidorDUKE},\footnote{Voisin's and Teixidor's cases complete each other quite remarkably.}
\cite{SchreyerLNM}, for curves with small Brill--Noether loci \cite{AproduMRL},
for plane curves \cite{LooseMANUSCRIPTA},
for curves on $K3$ surfaces \cite{VoisinJEMS,VoisinCOMP,Aprodu-FarkasCOMP},~etc.;
see also \cite{Aprodu-NagelULECT}
for a discussion.

\medskip

We shall consider in the sequel the case of curves on
$K3$ surfaces with emphasis on Voisin's approach to the problem
and the role played by Lazarsfeld--Mukai bundles.
It is interesting to notice that Green's conjecture for $K3$ sections can be formulated directly
in the $K3$ setup, as a vanishing result on the moduli space
$\mathcal{F}_g$ of polarized $K3$ surfaces. However, in the proof of this statement,
as it usually happens in mathematics, we have to exit the $K3$ world, prove
a more general result in the extended setup, and return to $K3$ surfaces.
The steps we have to take, ordered logically and not chronologically, are
the following.
In the first, most elaborated step, one finds an example for
odd genus \cite{VoisinCOMP,VoisinJEMS}. At this stage,
we are placed in the moduli space $\mathcal{F}_{2k+1}$.
Secondly, we exit the $K3$ world, land in $\mathcal{M}_{2k+1}$,
and prove the equality of two divisors \cite{Hirschowitz-RamananAENS,VoisinJEMS}.
The first step is used, and the identification of~the divisors extends
to their closure over the component $\varDelta_0$
of the boundary \cite{AproduMRL}.
In~the third step, we jump from a gonality stratum
$\mathcal{M}^1_{g,d}$ in a moduli space $\mathcal{M}_g$
to the~boundary of another moduli space of stable curves $\overline{\mathcal{M}}_{2k+1}$,
where $k=g-d+1$ \cite{AproduMRL}. The second step reflects into a vanishing result on
an explicit open subset of $\mathcal{M}^1_{g,d}$. Finally one goes back to
$K3$ surfaces and applies
the latter vanishing result \cite{Aprodu-FarkasCOMP} on $\mathcal{F}_g$.
In the steps
concerned with $K3$ surfaces (first and last), the Lazarsfeld--Mukai bundles are central objects.

\subsection{Voisin's Approach}
\label{subsection: Hilbert}

The proof of the generic Green conjecture was achieved by Voisin
in two papers \cite{VoisinJEMS,VoisinCOMP}, using a completely different
approach to Koszul cohomology via Hilbert scheme of points.

Let $X$ be
a  complex connected projective manifold and $L$ a line
bundle on $X$. It is obvious that any global
section $\sigma$ is uniquely determined by the
collection $\{\sigma(x)\}_x$, where $\sigma(x)\in L|_x\cong\mathbb{C}$
and $x$ belongs to a nonempty open subset of $X$.
One tries  to find a similar fact for
multisections in $\wedge^nH^0(X,L)$.

Let $\sigma_1\wedge\cdots\wedge\sigma_n$ be
a decomposable element in $\wedge^nH^0(X,L)$ with $n\ge 1$.
By analogy with the
case $n=1$, we have to look at the restriction
$\sigma_1|_{\xi}\wedge\cdots\wedge\sigma_n|_{\xi}\in \wedge^nL|_{\xi}$ where
$\xi$ is now a zero-dimensional subscheme, and it is clear
that we need $n$ points for otherwise this restriction would be zero.
Note that a zero-dimensional subscheme of length $n$ defines
a point in the punctual Hilbert scheme $X^{[n]}$. For technical
reasons, we shall restrict to curvilinear subschemes\footnote{A curvilinear
subscheme is defined locally, in the classical topology,  by $x_1=\cdots=x_{s-1}=x_s^k=0$;
equivalently, it is locally embedded in a smooth curve.}
which form a large open subset $X^{[n]}_c$ in a connected component
of the Hilbert scheme.\footnote{The connectedness of $X^{[n]}_c$
follows from the observation that a curvilinear subscheme is a deformation
of a reduced subscheme.}
Varying $\xi\in X^{[n]}_c$, the collection
$\{\sigma_1|_{\xi}\wedge\cdots\wedge\sigma_n|_{\xi}\}_{\xi}$
represents a section in a line bundle described as follows.
Put $\varXi_n\subset X^{[n]}_c\times X$ the incidence
variety and denote by $q$ and $p$ the projections on
the two factors; note that $q$ is finite of degree $n$.
Then $L^{[n]}:=q_*p^*(L)$ is a vector bundle of rank $n$
on $X^{[n]}_c$, and the fibre at a point $\xi\in X^{[n]}_c$
is $L^{[n]}|_{\xi}\cong L|_{\xi}$. In conclusion,
the collection $\{\sigma|_{\xi}\wedge\cdots\wedge\sigma|_{\xi}\}_{\xi}$
defines a section in the line bundle
$\mathrm{det}(L^{[n]})$.
The map we are looking at $\wedge^nH^0(L)\to H^0(\mathrm{det}(L^{[n]}))$
is deduced from the evaluation map
$\mathrm{ev}_n:H^0(L)\otimes \mathcal{O}_{X^{[n]}_c}\to L^{[n]}$,
taking $\wedge^n\mathrm{ev}_n$ and applying $H^0$.
It is remarkable that \cite{VoisinJEMS,VoisinCOMP,Ellingsrud-Goettsche-LehnJAG}:
\begin{thm}[Voisin, Ellingsrud--G\"ottsche--Lehn]
 The map\vspace*{-3pt}
 \[
 H^0(\wedge^n\mathrm{ev}_n):\wedge^nH^0(X,L)\to
 H^0\left(X^{[n]}_c,\mathrm{det}(L^{[n]})\right)\vspace*{-3pt}
 \]
 is an isomorphism.
\end{thm}

Since the exterior powers of $H^0(L)$ are building blocks
for Koszul cohomology, it is natural to believe that the
isomorphism above yields a relation between
the Koszul cohomology and the Hilbert scheme. To this
end, the Koszul differentials must  be reinterpreted in
the new context.

There is a natural birational morphism\footnote{We see
one advantage of working on $X^{[n]}_c$: subtraction makes sense only for curvilinear subschemes.}\vspace*{-3pt}
\[
\tau:\varXi_{n+1}\to X^{[n]}_c\times X,\
(\xi,x)\mapsto(\xi-x,x)\vspace*{-3pt}
\]
presenting $\varXi_{n+1}$ as the blowup of $X^{[n]}_c\times X$ along $\varXi_n$.
If we denote by $D_{\tau}$ the exceptional locus, we obtain
an inclusion \cite{VoisinJEMS}
\[
 q^*\mathrm{det}(L^{[n+1]})\cong\tau^*(\mathrm{det}(L^{[n]})\boxtimes
 L)(-D_{\tau})\hookrightarrow \tau^*(\mathrm{det}(L^{[n]})\boxtimes
 L)
\]
whence
\[
H^0\left(X^{[n+1]}_c,\mathrm{det}(L^{[n+1]})\right)\hookrightarrow
H^0(X^{[n]}_c\times X,\mathrm{det}(L^{[n]})\boxtimes L),
\]
identifying the left-hand member with the kernel of a Koszul
differential \cite{VoisinJEMS}. A~version of this identification
leads us to \cite{VoisinJEMS,VoisinCOMP}:

\begin{thm}[Voisin]
\label{thm: Hilbert scheme}
 For any integers $m$ and $n$, $K_{n,m}(X,L)$ is isomorphic to the
 cokernel of the restriction map:
 \[
 H^0\left(X^{[n+1]}_c\times X,\mathrm{det}(L^{[n+1]})\boxtimes L^{m-1}\right)\to
 H^0\left(\varXi_{n+1},\mathrm{det}(L^{[n+1]})\boxtimes L^{m-1}|_{\varXi_{n+1}}\right).
 \]
\end{thm}

The vanishing of Koszul cohomology is thus reduced  to proving
surjectivity of the restriction map above. In general, it
is very hard to prove surjectivity directly, and one
has to make a suitable base-change \cite{VoisinJEMS}.

\subsection{The Role of Lazarsfeld--Mukai Bundles in the Generic
Green Conjecture and Consequences}
\label{subsection: Role of LM}

In order to prove Green's conjecture for general curves,
it suffices to exhibit one example of a curve of
maximal Clifford index, which verifies the predicted
vanishing. Afterwards, the vanishing of
Koszul cohomology propagates by semicontinuity.
Even so, finding one
single example is a task of major difficulty.
The curves used by Voisin in \cite{VoisinJEMS,VoisinCOMP}
are $K3$ sections, and the setups change slightly, according to the parity
of the genus. For even genus, we have \cite{VoisinJEMS}:

\begin{thm}[Voisin]
 \label{thm: Voisin even}
Suppose that $g=2k$. Consider $S$ a $K3$ surface
with $\mathrm{Pic}(S)\cong \mathbb{Z}\cdot L$, $L^2=2g-2$, and
$C\in|L|$ a smooth curve. Then $K_{k,1}(C,K_C)=0$.
\end{thm}

For odd genus, the result is \cite{VoisinCOMP}:

\begin{thm}[Voisin]
 \label{thm: Voisin odd}
Suppose that $g=2k+1$. Consider $S$ a $K3$ surface
with $\mathrm{Pic}(S)\cong \mathbb{Z}\cdot L\oplus \mathbb{Z}\cdot \varGamma$, $L^2=2g-2$,
$\varGamma$ a smooth rational curve. $L\cdot \varGamma =2$ and
$C\in|L|$ a smooth curve. Then $K_{k,1}(C,K_C)=0$.
\end{thm}

Note that  the generic value for the Clifford index
in genus $g$ is $[(g-1)/2]$, and hence, in both cases, the prediction made by Green's
conjecture for general curve $C$ is precisely $K_{k,1}(C,K_C)=0$.

There are several reasons for making these choices:
the curves have maximal Clifford index, by
Theorem~\ref{thm: GL Cliff} (and the Clifford
dimension is one), the Lazarsfeld--Mukai
bundles associated to minimal pencils are $L$-stable,
the hyperplane section theorem applies, etc.

\medskip

We outline here the role played by Lazarsfeld--Mukai bundles in
Voisin's proof and, for simplicity, we restrict to the even-genus case.
By the \hbox{hyperplane} section Theorem~\ref{thm: Lefschetz}, the required
vanishing on the curve is equivalent to \hbox{$K_{k,1}(S,L)=0$.}
From the description of Koszul cohomology in terms of Hilbert
schemes, Theorem \ref{thm: Hilbert scheme}, adapting the notation from the previous subsection, one has to prove the surjectivity of the map
\[
q^*: H^0\left(S_c^{[n+1]},\mathrm{det}(L^{[n+1]})\right)\to
 H^0\left(\varXi_{n+1},q^*\mathrm{det}(L^{[n+1]})|_{\varXi_{n+1}}\right).
 \]
The surjectivity is proved after performing a suitable base-change.

We are in the case $\rho(g,1,k+1)=0$; hence there is a unique
Lazarsfeld--Mukai bundle $E$ on $S$ associated to all $g^1_{k+1}$
on curves in $|L|$. The uniqueness yields an alternate description
of $E$ as extension
\[
0\to \mathcal{O}_S\to E\to L\otimes I_{\xi}\to 0,
\]
where $\xi$ varies in $S_c^{[k+1]}$.

There exists a morphism
$\mathbb{P} H^0(S,E)\rightarrow S^{[k+1]}$ that sends a global section $s\in H^0(S,E)$
to its zero set $Z(s)$. By restriction to an open subset
$\mathbb{P}\subset\mathbb{P}
H^0(S,E)$, we obtain a morphism $\mathbb{P}\rightarrow S_c^{[k+1]}$,
inducing a commutative diagram
\begin{center}
                          \mbox{\xymatrix{
\mathbb{P}^{\prime} =  \mathbb{P}\times_{S_c^{[k+1]}}\varXi_{k+1}\ar[r]\ar[d]^{q^{\prime}} &\varXi_{k+1} \ar[d]^q \\
\mathbb{P} \ar[r] & S_c^{[k+1]}.
}}
                          \end{center}
Set-theoretically
\[
\mathbb{P}^{\prime} = \{(Z(s),x)|s\in H^0(S,E),x\in Z(s)\}.
\]

Unfortunately, this very natural base-change does not satisfy the
necessary conditions that imply the surjectivity
of $q^*$, \cite{VoisinJEMS}. Voisin modifies slightly this construction
and replaces $\mathbb{P}$ with another variety related to $\mathbb{P}$ which
parameterizes zero-cycles of the form
$Z(s)-x+y$
with $[s]\in\mathbb{P}$, $x\in\mathrm{Supp}(Z(s))$ and $y\in S$.
It turns out, after numerous elaborated calculations using the
rich geometric framework provided by the Lazarsfeld--Mukai
bundle, that the
new base-change is suitable and the surjectivity of $q^*$ follows
from vanishing results on the Grassmannian \cite{VoisinJEMS}.

\medskip

In the odd-genus case, Voisin proves first Green's conjecture
for smooth curves in $|L+\varGamma|$, which are easily seen to
be of maximal Clifford index. The situation on $|L+\varGamma|$
is somewhat close to the setup of Theorem~\ref{thm: Voisin even},
and the proof is similar. The next hard part is to descend
from the vanishing of $K_{k+1,1}(S,L\otimes\mathcal{O}_S(\varGamma))$
to the vanishing of $K_{k,1}(S,L)$. This step uses again
intensively the unique Lazarsfeld--Mukai bundle associated
to any $g^1_{k+2}$ on curves in $|L+\varGamma|$.

The odd-genus case is of maximal interest:
mixed with Hirschowitz-Ramanan result
\cite{Hirschowitz-RamananAENS}, Theorem~\ref{thm: Voisin odd}
gives a solution to Green's conjecture
for \textit{any} curve of odd genus and maximal Clifford index:

\begin{thm}[Hirschowitz--Ramanan, Voisin]
\label{thm: HRV}
 Let $C$ be a smooth curve of odd genus $2k+1\ge 5$
 and Clifford index $k$. Then $K_{k,1}(C,K_C)=0$.
\end{thm}

Note that Theorem~\ref{thm: HRV} implies
the following statement:

\begin{cor}
A smooth curve of odd genus
and maximal Clifford index
has Clifford dimension one.
\end{cor}

The proof of Theorem~\ref{thm: HRV}
relies on the comparison of two effective
divisors on the moduli space of curves $\mathcal{M}_{2k+1}$,
one given by the condition $\mathrm{gon}(C)\le k+1$, which
is known to be a divisor from \cite{Harris-MumfordINVENTIONES},
and the second given by $K_{k,1}(C,K_C)\ne 0$.
By duality $K_{k,1}(C,K_C)\cong K_{k-2,2}(C,K_C)$.
Note that $K_{k-2,2}(C,K_C)$ is isomorphic~to\vspace*{-3pt}
\[
\mathrm{Coker}
\left(\wedge^kH^0(K_C)\otimes H^0(K_C)/\wedge^{k+1}H^0(K_C)
\to H^0(\wedge^{k-1}M_{K_C}\otimes K_C^2)\right)\vspace*{-3pt}
\]
and the two members have the same dimension. The locus of curves
with $K_{k,1}\ne 0$ can be described as the degeneracy locus
of a morphism between vector bundles of the same dimension, and hence it
is a virtual divisor. Theorem~\ref{thm: Voisin odd} implies that this locus
is not the whole space, and in conclusion it must be an
effective divisor. \hbox{Theorem \ref{thm: GL nonvan}}
already gives an inclusion between the supports of
two divisors in question, and the
set-theoretic equality is obtained from a divisor class
calculation~\cite{Hirschowitz-RamananAENS}.

\vspace*{-6pt}
\subsection{Green's Conjecture for Curves on $K3$ Surfaces}
\label{subsection: Green for K3}

We have already seen that general $K3$ sections have
a mild behavior from the Brill--Noether theory viewpoint. In some
sense, they behave like general curves in any
gonality stratum of the moduli space of curves.

As in the previous subsections, fix a $K3$ surface $S$ and a globally
generated line bundle $L$ with $L^2=2g-2$ on $S$, and denote by $k$  the
gonality of a general smooth curve in the linear system $|L|$.
Suppose that $\rho(g, 1, k)\le 0$ to exclude the case $g=2k-3$
(when $\rho(g, 1, k)=1$).
If in addition the curves in $|L|$ have Clifford dimension one,
Theorem~\ref{thm: Green Cliffdim 1}
shows that\vspace*{-3pt}
\[
\mathrm{dim}\ W^1_{g-k+2}(C)=\rho(g,1,g-k+2)=g-2k+2,\vspace*{-3pt}
\]
property which was called the \textit{linear growth condition}.

This property appears in connection with Green's conjecture~\cite{AproduMRL}
for a much larger class of curves:

\begin{thm}
\label{thm: Aprodu}
 If $C$ is any smooth curve of genus $g\ge 6$ and gonality $3\le k<[g/2]+2$
 with $\mathrm{dim}\ W^1_{g-k+2}(C)=\rho(g,1,g-k+2)$,
 then $K_{g-k+1,1}(C,K_C)=0$.
\end{thm}

One effect of Theorems~\ref{thm: Aprodu} and~\ref{thm: GL nonvan} is that an arbitrary curve that satisfies the
linear growth condition is automatically of Clifford
dimension one and verifies Green's conjecture.

Theorem~\ref{thm: Aprodu} is a consequence of Theorem~\ref{thm: HRV} extended over the boundary of the moduli space.
Starting from a $k$-gonal smooth curve $[C]\in \mathcal{M}_g$, by identifying pairs of general points
$\{x_i,y_i\}\subset C$ for $i=0, \dots, g-2k+2$ we produce a stable
irreducible curve
\[
\left[X:=C/(x_0\sim y_0, \ldots,   x_{g-2k+2}\sim y_{g-2k+2})\right]
\in  \overline{\mathcal{M}}_{2(g-k+1)+1},
\]
and the Koszul cohomology of $C$ and of $X$
are related by the inclusion $K_{p,1}(C,K_C)$ $\subset K_{p,1}(X,\omega_X)$
for all $p\ge 1$, \cite{VoisinJEMS}.
If $C$ satisfies the linear growth condition then $X$
has maximal gonality\footnote{The gonality for
a singular stable curve is defined in terms of admissible
covers \cite{Harris-MumfordINVENTIONES}.} $\mathrm{gon}(X)=g-k+3$, i.e., $X$ lies outside the
closure of the divisor $\mathcal{M}_{2(g-k+1)+1, g-k+2}^1$
consisting of curves with a pencil
$g^1_{g-k+2}$. The class of the failure locus of
Green's conjecture on $\overline{\mathcal{M}}_{2(g-k+1)+1}$
is a multiple of the divisor $\overline{\mathcal{M}}_{2(g-k+1)+1, g-k+2}^1$;
hence Theorem~\ref{thm: HRV} extends to irreducible stable curves
of genus $2(g-k+1)+1$ of maximal gonality $(g-k+3)$.
In particular, $K_{g-k+1,1}(X,\omega_X)=0$, implying
$K_{g-k+1,1}(C,K_C)=0$.

\medskip

Coming back to the original situation, we conclude from
Theorems~\ref{thm: Aprodu} and~\ref{thm: Green Cliffdim 1}
and Corollary~\ref{cor: Lefschetz K3}
that Green's conjecture holds for a $K3$ section $C$
having Clifford dimension one.
If $\mbox{Cliff}(C)=\mathrm{gon}(C)-3$,
either $C$ is a smooth plane curve or else there
exist smooth curves $D, \varGamma \subset S$, with
$\varGamma^2=-2, \varGamma\cdot D=1$ and $D^2\geq 2$, such that $C\equiv 2D+
\varGamma$ and $\mbox{Cliff}(C)=\mbox{Cliff}(\mathcal{O}_C(D))$
\cite{Ciliberto-PareschiCRELLE,KnutsenIJM}. The
linear growth condition is no longer satisfied,
and this case is treated differently, by
degeneration to a reduced curve
with two irreducible components \cite{Aprodu-FarkasCOMP}.

\medskip

The outcome of this analysis of the Brill--Noether loci
is the following \cite{VoisinJEMS,VoisinCOMP,Aprodu-FarkasCOMP}:

\begin{thm}
\label{thm: Green on K3}
Green's conjecture is valid for any smooth curve on a $K3$ surface.
\end{thm}

Applying Theorem~\ref{thm: Green on K3},
Theorem~\ref{thm: Lefschetz}, and the duality,
we obtain a full description of the situations when Koszul
cohomology of a $K3$ surface is zero \cite{Aprodu-FarkasCOMP}:

\begin{thm}
\label{thm: K3}
Let $S$ be a $K3$ surface and $L$ a globally generated line bundle
with $L^2=2g-2\ge 2$. The Koszul cohomology group $K_{p,q}(S,L)$ is nonzero if and only if one of the following cases occurs:
\begin{itemize}
\item[(1)] $q=0$ and $p=0$, or
\item[(2)] $q=1$, $1\le p\le g-c-2$, or
\item[(3)] $q=2$ and $c\le p\le g-1$, or
\item[(4)] $q=3$ and $p=g-2$.
\end{itemize}
\end{thm}

The moral is that the shape of the Betti table, i.e., the distribution of
zeros in the table, of a polarized $K3$ surface is completely
determined by the geometry of hyperplane sections;
this is one of the many situations where algebra and geometry are intricately related.

\section{Counterexamples to Mercat's Conjecture in Rank Two}
\label{section: Higher BN}

Starting from Mukai's works,  experts tried
to generalize the classical Brill--Noether
theory to higher-rank vector bundles on curves.
Within these extended theories,\footnote{Higher-rank
Brill--Noether theory is a major, rapidly growing research field,
and it deserves a separate dedicated survey.}
we note the attempt
to find a proper generalization of the Clifford index.
H. Lange and P. Newstead proposed the following
definition. Let $E$ be a semistable vector
bundle of rank $n$ of degree $d$ on a smooth curve $C$.
Put
\[
 \gamma(E):=\mu(E)-2\frac{h^0(E)}{n}+2.
\]

\begin{defn}[Lange--Newstead]
 The \textit{Clifford index} of rank $n$ of $C$ is
\[
 \mathrm{Cliff}_n(C):=\mathrm{min}\{\gamma(E):\ \mu(E)\le g-1,\ h^0(E)\ge 2n\}.
\]
\end{defn}

From the definition, it is clear that $\mathrm{Cliff}_1(C)=\mathrm{Cliff}(C)$
and $\mathrm{Cliff}_n(C)\le\mathrm{Cliff}(C)$ for all $n$.\footnote{For
any line bundle $A$, we have $\gamma(A^{\oplus n})=\mathrm{Cliff}(A)$.}

Mercat conjectured \cite{MercatIJM} that $\mathrm{Cliff}_n(C)=\mathrm{Cliff}(C)$.
In rank two, the conjecture is known to hold in a number of cases: for general curves of
small gonality, i.e., corresponding to a general point in a gonality
stratum $\mathcal{M}^1_{g,k}$ for small $k$ (Lange-Newstead), for plane curves
(Lange--Newstead), for general curves of genus $\le 16$
(Farkas--Ortega), etc.
However, even in rank two, the conjecture is false. It is remarkable
that counterexamples are found for curves of maximal Clifford index~\cite{Farkas-OrtegaIJM}:

\begin{thm}[Farkas--Ortega]
 Fix $p\ge 1$, $a\ge 2p+3$. Then there exists a smooth curve
of genus $2a+1$ of maximal Clifford index lying on a smooth $K3$
surface $S$ with $\mathrm{Pic}(S)=\mathbb{Z}\cdot C\oplus \mathbb{Z}\cdot H$,
$H^2=2p+2$, $C^2=2g-2$, $C\cdot H=2a+2p+1$, and there exists a
stable rank-two vector bundle $E$ with $\mathrm{det}(E)=\mathcal{O}_S(H)$
with $h^0(E)=p+3$, $\gamma(E)=a-\frac{1}{2}<a=\mathrm{Cliff}(A)$, and
hence Mercat's conjecture in rank two fails for~$C$.
\end{thm}

The proof uses restriction of Lazarsfeld--Mukai bundles.
However, it is interesting that the bundles are not
restricted to the same curves to which they
are associated. More precisely, the genus of
$H$ is $2p+2$ and $H$ has maximal gonality $p+2$.
Consider $A$ a minimal pencil on $H$, and take
$E=E_{H,A}$ the associated Lazarsfeld--Mukai bundle.
The restriction of $E$ to $C$ is stable and verifies
all the required properties.

\medskip

A particularly interesting case is $g=11$.
In this case, as shown by Mukai \cite{MukaiLNPAM},
a general curve $C$ lies on a unique $K3$
surface $S$ such that $C$ generates $\mathrm{Pic}(S)$.\vadjust{\pagebreak}
It~is~remarkable that the failure locus of Mercat's conjecture in rank
two \textit{coincides} with the Noether-Lefschetz divisor
\[
\mathcal{NL}^4_{11,13}:=\left\{[C]\in \mathcal{M}_{11}:
\begin{array}{l}
C\mbox{ lies on a } K3 \mbox{ surface } S,  \
\mathrm{Pic}(S)\supset \mathbb{Z}\cdot C\oplus \mathbb{Z}\cdot H,\\
H\in \mathrm{Pic}(S) \mbox{ is nef},
 H^2=6, \ C\cdot H=13, \ C^2=20 \end{array}
\right\}
\]
inside the moduli space $\mathcal{M}_{11}$.
We refer to \cite{Farkas-OrtegaIJM} for details.


\begin{thebibliography}{99}
 \bibitem{AproduMRL}
 Aprodu, M.:
 Remarks on syzygies of d-gonal curves.
 Math. Res. Lett. {\bf 12}(3), 387--400 (2005)

  \bibitem{Aprodu-NagelULECT}
 Aprodu, M., Nagel, J.:
 Koszul cohomology and algebraic geometry.
 University Lecture Series, vol. 52. American Mathematical Society, Providence (2010)

 \bibitem{Aprodu-FarkasCOMP}
 Aprodu, M., Farkas, G.:
 Green's conjecture for curves on arbitrary K3 surfaces.
 Compositio Math. {\bf 147}, 839--851  (2011)


 \bibitem{Arbarello-Cornalba-Griffiths-Harris}
 Arbarello, E., Cornalba, M., Griffiths, P.A., Harris, J.:
 Geometry of algebraic curves, vol. I.
 Grundlehren der Mathematischen Wissenschaften, vol. 267.
 Springer, New York (1985)

 \bibitem{Ciliberto-PareschiCRELLE}
 Ciliberto, C., Pareschi, G.:
 Pencils of minimal degree on curves on a $K3$ surface.
 J. Reine Angew. Math. {\bf 460}, 15--36  (1995)

 \bibitem{Donagi-MorrisonJDG}
 Donagi, R., Morrison, D.R.:
 Linear systems on $K3$ sections. Diff. J. Geom. {\bf 29} 49--64  (1989)

 \bibitem{EisenbudBOOK}
  Eisenbud, D.: Geometry of syzygies.
 Graduate Texts in Mathematics, vol. 229. Springer, New York (2005)

 \bibitem{ELMS}
 Eisenbud, D., Lange, H., Martens, G.,  Schreyer, F.-O.:
 The Clifford dimension of a projective curve.
 Compositio Math. {\bf 72}, 173--204  (1989)

 \bibitem{Ellingsrud-Goettsche-LehnJAG}
 Ellingsrud, G., G\"ottsche, L., Lehn, M.:
 On the cobordism class of the Hilbert scheme of a surface.
 Alg. J.: Geom. {\bf 10},  81--100  (2001)

 \bibitem{Farkas-OrtegaIJM}
 Farkas, G., Ortega, A.:
 Higher-rank Brill--Noether theory on sections of $K3$ surfaces.
 Internat. Math. J. {\bf 23}, 1250075, 18 pp (2012)

 \bibitem{GreenJDG}
 Green, M.:
 Koszul cohomology and the geometry of projective varieties.
 Diff. J. Geom. {\bf 19}, 125--171  (1984)

  \bibitem{Green-LazarsfeldINVENTIONES}
 Green, M., Lazarsfeld, R.:
 Special divisors on curves on a $K3$ surface.
 Inventiones Math. {\bf 89}, 73--90  (1987)

 \bibitem{GreenICTP}
 Green, M.:  Koszul cohomology and geometry.
 In: Cornalba, M., et al. (eds.) Proceedings of the First College on Riemann Surfaces Held in Trieste, Italy,
 November 1987,  pp. 177--200. World Scientific, Singapore (1989)

 \bibitem{Harris-MumfordINVENTIONES}
 Harris, J., Mumford, D.:
 On the Kodaira dimension of the moduli space of curves.
 Inventiones Math. {\bf 67},  23--86   (1982)

 \bibitem{Hirschowitz-RamananAENS}
 Hirschowitz, A., Ramanan, S.:
 New evidence for Green’s conjecture on syzygies of canonical curves.
 Ann. Sci. \'Ecole Norm. Sup. {\bf 31}, 145--152  (1998)

 \bibitem{Iskovskih-Prokhorov}
 Iskovskih, V.A., Prokhorov, Yu.:
 Fano varieties. In:  Parshin, A.N., Shafarevich, I.R. (eds.) Algebraic Geometry V. Encyclopedia of Mathematical Science.
 Springer, New York (1999)

 \bibitem{KnutsenIJM}
 Knutsen, A.:
 On two conjectures for curves on $K3$ surfaces.
 Internat. Math. J. {\bf 20}, 1547--1560  (2009)

 \bibitem{LazarsfeldJDG}
 Lazarsfeld, R.:
 Brill--Noether--Petri without degenerations.
 Diff. J. Geom. {\bf 23}, 299--307  (1986)

 \bibitem{LazarsfeldICTP}
 Lazarsfeld, R.:
 A sampling of vector bundle techniques in the study of linear series.
 In: Cornalba, M., et al. (eds.) Proceedings of the first college on Riemann surfaces held in Trieste, Italy,
 November 1987, pp. 500--559. World Scientific, Singapore (1989)

 \bibitem{Lelli-Chiesa}
 Lelli-Chiesa, M.:
 Stability of rank-3 Lazarsfeld--Mukai bundles on K3 surfaces.
 Proc. Lond. Math. Soc. {\bf 107}, 451--479 (2013)
 
 \bibitem{LooseMANUSCRIPTA}
 Loose, F.:
 On the graded Betti numbers of plane algebraic curves.
 Manuscr. Math. {\bf 64}, 503--514   (1989)


 \bibitem{MercatIJM}
 Mercat, V.:
 Clifford's theorem and higher rank vector bundles.
 Internat. Math. J.: {\bf 13}, 785--796 (2002)

 \bibitem{MukaiPNAS}
 Mukai, S.:
 Biregular classification of Fano 3-folds and Fano manifolds of coindex 3.
 Proc. Natl. Acad. Sci. USA {\bf 86}, 3000--3002  (1989)

 \bibitem{MukaiLMS}
 Mukai, S.:
 Fano 3-folds.
 London Math. Soc. Lect. Note Ser. {\bf 179}, 255--263  (1992)

 \bibitem{MukaiLNPAM}
 Mukai, S.:  Curves and $K3$ surfaces of genus eleven.
 In: Moduli of Vector Bundles. Lecture Notes in Pure Applied Mathematics,
 vol. 179, pp. 189--197. Dekker, New York (1996)

 \bibitem{PareschiJAG}
 Pareschi, G.:
 A proof of Lazarsfeld's theorem on curves on $K3$ surfaces.
 Alg. J. Geom. {\bf 4}, 195--200 (1995)

 \bibitem{ReidJLM}
 Reid, M.:
 Special linear systems on curves lying on $K3$ surfaces.
 J. London Math. Soc. (2) {\bf 13}, 454--458  (1976)

 \bibitem{Saint-DonatAJM}
 Saint-Donat, B.:
 Projective models of $K3$ surfaces.
 American J. Math. {\bf 96}(4), 602--639  (1974)

 \bibitem{SchreyerLNM}
 Schreyer, F.-O.:
 Green's conjecture for general $p$-gonal curves of large genus.
 Algebraic curves and projective geometry (Trento, 1988),
 pp. 254--260. Lecture Notes in Mathematics, vol. 1389. Springer,
 Berlin (1989)

 \bibitem{TeixidorDUKE}
 Teixidor i Bigas, M.:
 Green's conjecture for the generic $r$-gonal curve of genus $g\geq 3r-7$.
 Duke Math. J. {\bf 111}, 195--222  (2002)

 \bibitem{VoisinJEMS}
 Voisin, C.:
 Green's generic syzygy conjecture for curves of even genus lying on a $K3$ surface.
 J.  European Math. Soc. {\bf 4}, 363--404  (2002)

 \bibitem{VoisinCOMP}
 Voisin, C.:
 Green's canonical syzygy conjecture for generic curves of odd genus.
 Compositio Math. {\bf 141}, 1163--1190  (2005)

\end{thebibliography}
\end{document}